\documentclass[a4paper,10pt]{article}
\usepackage[cp1251]{inputenc}
\usepackage{enumerate}
\usepackage{indentfirst}
\usepackage[font = scriptsize, position = below, aboveskip = 5pt]{caption}
\usepackage{subfigure}
\usepackage{graphicx}
\usepackage{algorithmic}
\usepackage{amsmath, amssymb}
\usepackage{color}
\usepackage{wrapfig}

\newtheorem{remark}{Remark}
\newtheorem{lemma}{Lemma}

\newtheorem{theorem}{Theorem}

\newenvironment{Proof} 
{\par\noindent{\bf Proof.}\newline} 
{\hfill$\scriptstyle\blacksquare$\newline}

\title{About new dynamical interpretations of entropic model of 
correspondence matrix calculation and Nash-Wardrop's equilibrium in 
Beckmann's traffic flow distribution model}
\author{Gasnikova E.\,V.\and
Nagapetyan T.\.A.}

\begin{document}
\maketitle

\section[Correspondence Matrix Calculation Model] {Correspondence Matrix Calculation Model\footnote{This section was written, based on the works \cite{1} -- \cite{11}.}}
Assume that in some town there are $n$ districts (regions), $L_i >0$ is the number of 
residents living at the district $i$, and $W_j >0$ is the number of residents working at the district $j$. 
By $x_{ij} \left( t \right)\ge 0$ we will denote number of residents living at the district $i$ and working at the district $j$ 
at the moment of time $t\ge 0$. Over course of the time numbered residents (whose quantity doesn't 
alter, and is equal to $N=\sum\limits_{i=1}^n {L_i } =\sum\limits_{j=1}^n {W_j }$) change their apartments (homes). And we suppsose, that changes may happen 
only by the means of the exchange of the apartments, i.e.
\begin{equation}
\label{A}
x_{ij} \left( t \right)\ge 0,
\quad
\sum\limits_{j=1}^n {x_{ij} \left( t \right)} \equiv L_i ,
\quad
\sum\limits_{i=1}^n {x_{ij} \left( t \right)} \equiv W_j ,
\quad
i,j=1,\ldots,n. 
\end{equation}

Suppose, that at the time $t\ge 0$ resident $r$ lives in the district $k$ and 
works in the district $m$, and the resident $l$ lives in the district $p$ and works in 
the district $q$. Then $p_{k,m;\;p,q}^L \left( t \right)\Delta t+o\left( {\Delta 
t} \right)$ - is the probability for the residents with numbers $r$ and $l$ $(1\le 
r<l\le N)$ to exchange their apartments in a period of time $\left({t,t+\Delta t} \right)$. 
It's natural to consider that probability (in the unit time) of exchanging apartments depends only on the location of the workplaces and 
homes, which are exchanged. For instance, it may be considered that the 
``distance'' between district $i$ and district $j$ is $c_{ij} \ge 0$, and
\[
p_{k,m;\;p,q}^L \left( t \right)\equiv p^L\exp \;(\underbrace {\left( 
{c_{km} +c_{pq} } \right)}_{\begin{array}{c}
 \mbox{\footnotesize{sum of the distances}} \\ 
 \mbox{\footnotesize{before exchange}} \\ 
 \end{array}}-\underbrace {\left( {c_{pm} +c_{kq} } 
\right)}_{\begin{array}{c}
 \mbox{\footnotesize{sum of the distances}} \\ 
 \mbox{\footnotesize{after exchange}} \\ 
 \end{array}}\;)>0.
\]
Then, by the virtue of the ergodic theorem for discrete homogeneous Markov process 
with finite number of states, for all $\left\{ {x_{ij} } \right\}_{i=1,j=1}^{n,\,n}\in$(\ref{A}) we have that 
$$\mathop {\lim }\limits_{t\to \infty } P\left( {x_{ij} \left( t \right)=x_{ij} ,i,j=1,\ldots,n} \right)=$$
$$=Z^{-1}\prod\limits_{i=1,j=1}^{n,\;n} {\exp \left( {-2c_{ij} x_{ij} } \right)\cdot \left( {x_{ij} !} \right)^{-1}} \mathop =\limits^{def} p\left( {\left\{ {x_{ij} } \right\}_{i=1,\;j=1}^{n,\;n} } \right),$$
where $Z$ is the normalizing multiplier.

Here we have the case where the final distribution, which is also a stationary distribution, satisfies detailed balance condition:\footnote{ Multipliers before 
probabilities, for example, in the state $\left\{ {x_{ij} } 
\right\}_{i=1,\;j=1}^{n,\;n} $, arise because of the number of the ways to 
choose the resident, living in the district $p$ and working in the district $m$, is 
$x_{pm} $, and independently the number of ways to choose the resident, living 
in the district $k$ and working in the district $q$, is $x_{kq} $.}
$$
\left( {x_{km} +1} \right)\left( {x_{pq} +1} \right)\cdot\hat{p}\cdot p_{k,m;\;p,q}^L =x_{pm} x_{kq} p\left( {\left\{ {x_{ij} } \right\}_{i=1,\;j=1}^{n,\;n} } 
\right)p_{p,m;\;k,q}^L, 
$$
where $$\hat{p}=p\left(\left\{ {x_{11} ,\ldots,x_{km} +1,\ldots,x_{pq} +1,\ldots,x_{pm} -1,\ldots,x_{kq} -1,\ldots,x_{nn} } \right\} \right).$$
Distribution $p\left( {\left\{ {x_{ij} } \right\}_{i=1,\;j=1}^{n,\;n} } 
\right)$ on a set (\ref{A}) is concentrated with $N\gg 1$ (see below) in a 
neighborhood of the most probable value $\left\{ {x_{ij}^\ast } 
\right\}_{i=1,\;j=1}^{n,\;n} $, which is determined as a solution of the following entropic --linear programming problem:

\begin{equation}
\label{EP}
 \sum\limits_{i=1,\;j=1}^{n,\;n} {x_{ij} \ln x_{ij} } +2\sum\limits_{i=1,\;j=1}^{n,\;n} {c_{ij} x_{ij} } \to \mathop {\min 
}\limits_{\left\{ {x_{ij} } \right\}_{i=1,\;j=1}^{n,\;n} \in \left(\ref{A} \right)} .
\end{equation}

Solution of this problem might be presented as 
$$x_{ij} =\exp \left( {-1-\lambda _i^L -\lambda _j^W -2c_{ij} } \right),$$ 
where Lagrange multipliers (dual variables) $\left\{ {\lambda _i^L } \right\}_{i=1}^n $ and 
$\left\{ {\lambda _j^W } \right\}_{j=1}^n $ are determined\footnote{ This can be done in a different ways. For example, 
by Bregman's balancing method or by Newton's method \cite{5}. The other way is to solve the dual problem for the entropy programming problem (\ref{EP}). 
There are a lot of different algorithms with the first order oracle (MART, GISM, etc. \cite{4}, \cite{5}). It can be shown that most of this methods (including Bregman's) are just 
barrier-multiplicative antigradient descending methods \cite{11}. At the end (when the iteration process is achieving a sufficient small vicinity of the global minimum) it is 
worth to use the second order interior-point method, like Nesterov--Nemirovskii polynomial algorithm \cite{12} (for so-called ``separable'' tasks).}  from the system of equations (\ref{A}).
In practice we usually have some information about 
$\left\{ {L_i ,W_i } \right\}_{i=1}^n $ and $\left\{ {c_{ij} } 
\right\}_{i=1,j=1}^{n,\;n} $. So, when we solve (\ref{EP}), we find 
$$x_{km} \left( {\left\{ {L_i ,W_i } \right\}_{i=1}^n ;\left\{ {c_{ij} } \right\}_{i=1,j=1}^{n,\;n} } \right)$$

If $N\sim nm$, $m\gg1$ $\forall i,j=1,\ldots,n\to L_i ,W_j \sim m,\,c_{ij}=c>0$, then the distribution of the probabilities 
$p\left( {\left\{ {x_{ij} } \right\}_{i=1,\;j=1}^{n,\;n} } \right)$ on the 
set (\ref{A}) is concentrated in ${\rm O}\left(\sqrt{m}\right)$ neighborhood of the most probable 
value $x_{ij}^\ast \approx {L_i W_j } \mathord{\left/ {\vphantom {{L_i W_j } 
N}} \right. \kern-\nulldelimiterspace} N\sim m/n,\,i,j=1,\ldots,n$. More 
precisely:
$$
\exists \lambda >0\colon\mathop {\lim }\limits_{t\to \infty } P\left(
\forall\,i,j=1,\ldots,n\to \left| {{x_{ij} \left( t \right)} \mathord{\left/ {\vphantom {{x_{ij} \left( t \right)} {x_{ij}^\ast }}} \right. \kern-\nulldelimiterspace} {x_{ij}^\ast }-1} \right|\le \lambda/\sqrt{m}\right)\ge 0.999
$$
Indeed, let us note, that
$$
\forall \;\;\left\{ {x_{ij} } \right\}_{i=1,\;j=1}^{n,\;n} \in (\ref{A})\to \sum\limits_{i=1,\;j=1}^{n,\;n} {\frac{\partial \ln 
p\left( {\left\{ {x_{ij}^\ast } \right\}_{i=1,\;j=1}^{n,\;n} } 
\right)}{\partial x_{ij} }\cdot \left( {x_{ij} -x_{ij}^\ast } \right)} \le 
0
$$
Thus, $\forall \;\;\left\{ {x_{ij} } \right\}_{i=1,\;j=1}^{n,\;n} \in (\ref{A})\;\;\exists \;\;\theta \in \left[ {0,1} \right]:$

\begin{multline*}
\ln p\left( {\left\{ {x_{ij} } \right\}_{i=1,\;j=1}^{n,\;n} } \right)\le \ln 
p\left( {\left\{ {x_{ij}^\ast } \right\}_{i=1,\;j=1}^{n,\;n} } 
\right)+\\
+\sum\limits_{i=1,\;j=1}^{n,\;n} {\frac{\partial ^2\ln p\left( 
{\left\{ {x_{ij}^\ast \theta +x_{ij} \cdot \left( {1-\theta } \right)} 
\right\}_{i=1,\;j=1}^{n,\;n} } \right)}{\partial x_{ij}^2 }\cdot 
\frac{\left( {x_{ij} -x_{ij}^\ast } \right)^2}{2}} 
\end{multline*}

Since
\[
\frac{\partial ^2\ln p\left( {\left\{ {x_{ij} } \right\}_{i=1,\;j=1}^{n,\;n} 
} \right)}{\partial x_{ij}^2 }=\frac{\partial ^2\left( 
{-\sum\limits_{i=1,\;j=1}^{n,\;n} {x_{ij} \ln \;x_{ij} } } \right)}{\partial 
x_{ij}^2 }=-\frac{1}{x_{ij} }
\]
we have ``inequality of measure concentration'':
\[
\forall \;\;M>0,
\quad
\forall \;\;\left\{ {x_{ij} } \right\}_{i=1,\;j=1}^{n,\;n} \in (\ref{A}):
\quad
\sum\limits_{i=1,\;j=1}^{n,\;n} {\frac{\left( {x_{ij} -x_{ij}^\ast } 
\right)^2}{2\max \left\{ {x_{ij} ,\;x_{ij}^\ast } \right\}}} \ge M
\]
\[
p\left( {\left\{ {x_{ij} } \right\}_{i=1,\;j=1}^{n,\;n} } \right)\le 
e^{-M}p\left( {\left\{ {x_{ij}^\ast } \right\}_{i=1,\;j=1}^{n,\;n} } 
\right)
\]
\section[Beckmann traffic flow distribution model]{Beckmann traffic flow distribution model\footnote{This section was written, based on the works \cite{13} -- \cite{24}.}}

Let us consider the oriented graph $\Gamma =\left( {V,E} \right)$, which stands for transportation 
route in some town ($V$ -- nodes (vertices), $E\subset V\times V$ -- arc of 
the network (edges)). Let $W=\left\{ {w=\left( {i,j} \right):\;i,j\in V} 
\right\}$ be a set of pairs inlet-outlet; $p=\left\{ {v_1 ,v_2 ,\ldots,v_m } 
\right\}$ -- route from $v_1 $ to $v_m $, if $\left( {v_k ,v_{k+1} } 
\right)\in E$, $k=1,\ldots,m-1$ (it will be shown later (see example by V.I. 
Shvetsov) that, to specify the path it may not be enough to indicate only 
the set of vertices. In general, one must also specify exactly which edge, 
connecting the specified vertices, is chosen); $P_w $ -- set of routes in 
correspondence $w\in W$; $P=\bigcup\nolimits_{w\in W} {P_w } $ -- collection 
of all routes in the network $\Gamma $; $x_p $ -- flux on the way $p$, $\vec{x}=\left\{ {x_p :\;\;p\in P} \right\}$; $G_p \left( {\vec{x}} \right)$ -- 
specific costs of travel on the road $p$, $\vec{G}\left( {\vec{x}} 
\right)=\left\{ {G_p \left( {\vec{x}} \right):\;\;p\in P} \right\}$; $y_e $ 
-- flux on the arc $e$: $\vec{y}=\Theta \vec{x}$, where $\Theta =\left\{ 
{\delta _{ep} } \right\}_{e\in E,p\in P} $ ($\delta _{ep} =\left\{ {1,\;e\in 
p;\;0,\;e\notin p} \right\})$; $\tau _e \left( {y_e } \right)$ -- specific 
costs of travel on the arc $e$ (generally increasing, convex, smooth 
functions), it is natural to assume, that $\vec{G}\left( {\vec{x}} 
\right)=\Theta ^T\vec{\tau }\left( {\vec{y}} \right)$. Let flows on 
correspondences $d_w $, $w\in W$ to be known. Then $\vec{x}$, which 
describes flow distribution, must lie in the set:
$$
X=\left\{ \vec{x}\ge 0\colon\sum\limits_{p\in P_w } x_p  =d_w ,\;w\in W \right\}.
$$
Consider a game in which each element $w\in W$ corresponds to a considerably 
big ($d_w \gg 1)$ set of players of the same type. The set of pure 
strategies of such player is $P_w $, and profit (minus losses) is defined by 
the formula $-G_p \left( {\vec{x}} \right)$ (a player chooses a strategy 
$p\in P_w $ and neglects the fact, that $\left| {P_w } \right|$ components 
of the vector $\vec{x}$ and hence the profit depends slightly on his 
choice). One can show, that Nash equilibrium is equivalent to 
complementarity problem, which equivalent to a solution of variation 
inequality, which, in its turn equivalent to a solution of convex 
optimization problem.
$$
\forall w\in W,\,p\in P_w \to x_p^\ast \,\cdot \left( {G_p \left( {\vec{x}^\ast } \right)-\mathop {\min }\limits_{q\in P_w } \;G_q \left( {\vec{x}^\ast } \right)} \right)=0
$$
$$\Updownarrow$$
$$
\forall \vec{x}\in X\to \left\langle {\vec{G}\left( {\vec{x}^\ast } 
\right),\vec{x}-\vec{x}^\ast } \right\rangle \ge 0
$$
$$\Updownarrow$$
\begin{equation}
 \label{eq:psi}
\Psi \left( {\vec{x}} \right)=\sum\limits_{e\in E} 
{\int\limits_0^{\sum\nolimits_{p\in P} {x_p \delta _{pe} } } {\tau _e \left( 
z \right)dz} } \to \mathop {\min }\limits_{\vec{x}\in X} .
\end{equation}
It is easy to show, that in the case $\vec{G}\left( {\vec{x}} \right)$ being 
strictly monotonic transformation, the Nash-Vardrop equilibrium $\vec{x}^\ast $ is unique. If $\tau _e \left( {y_e } \right)$ are increasing 
functions then $\vec{y}^\ast =\Theta \vec{x}^\ast $ is unique, although, as 
we will see later, $\vec{x}^\ast $ isn't necessarily unique.

The route at a step $\left( {n+1} \right)$ player at correspondence $w$, 
choose independently according the mixed strategy with probability
\[
\mbox{Prob}_p^w \left( {n+1} \right)={\gamma _n \cdot \max \left\{ {x_p 
\left( n \right),n^{-1}} \right\}\exp \left( {{-G_p \left( {\vec{x}\left( n 
\right)} \right)} \mathord{\left/ {\vphantom {{-G_p \left( {\vec{x}\left( n 
\right)} \right)} T}} \right. \kern-\nulldelimiterspace} T} \right)} 
\mathord{\left/ {\vphantom {{\gamma _n \cdot \max \left\{ {x_p \left( n 
\right),n^{-1}} \right\}\exp \left( {{-G_p \left( {\vec{x}\left( n \right)} 
\right)} \mathord{\left/ {\vphantom {{-G_p \left( {\vec{x}\left( n \right)} 
\right)} T}} \right. \kern-\nulldelimiterspace} T} \right)} {Z_n^w }}} 
\right. \kern-\nulldelimiterspace} {Z_n^w },\;\;w\in W,
\]
to choose path $p\in P_w $ ($0<\gamma _n \le 1)$, and with probability 
$1-\gamma _n $ to choose the same strategy as at the \textit{n-th }step. Here 
$x_p \left( n \right)$ -- number of players at $w$, who have chosen at the 
\textit{n-th }step strategy $p\in P_w $, and $Z_n^w $ can be found from the normalization 
condition. Multiplier $\max \left\{ {x_p \left( n \right),n^{-1}} \right\}$ 
describes the will to imitate and, also, the reliability of using this 
strategy. This multiplier notices specifics of the problem (without it there 
could be convergence to something different from the Nash-Wardrop equilibrium). Parameter $\gamma $ 
describes ``the conservatism'', while ``the temperature'' $T$ stands for ``the risk 
appetite''.

\begin{theorem}
\label{th:1}
Let $T>0$ be sufficiently small, $\sum\limits_{n=1}^\infty {\gamma _n } =\infty , \sum\limits_{n=1}^\infty {\left( {\gamma _n } \right)^2} <\infty $. Then $\Psi 
\left( {\vec{x}\left( n \right)} 
\right)\mathrel{\mathop{\kern0pt\longrightarrow}\limits_{n\to \infty 
}^{a.s.}} \Psi \left( {\vec{x}^\ast } \right)$, where $\vec{x}^\ast$ is the minimizer from \eqref{eq:psi}. Moreover, if the equilibrium is unique, then $\vec{x}\left( n 
\right)\mathrel{\mathop{\kern0pt\longrightarrow}\limits_{n\to \infty}^{a.s.}} \vec{x}^\ast  .$
\end{theorem}

In the experiments, conducted at the Laboratory of Experimental Economics in the Faculty of Applied Mathematics and Control, MIPT, in which students of the 5$^{th}$ course were 
involved, we observed the convergence to equilibrium and ``vibrations'' 
around it. Fluctuations should be explained, apparently, by the fact that in 
experiments the number of players was small and the hypothesis of a 
competitive market was not performed. We also observed, that for students 
$\gamma _n \equiv \gamma >0$ it is more likely, than $\gamma _n \sim 1 
\mathord{\left/ {\vphantom {1 n}} \right. \kern-\nulldelimiterspace} n$. As 
a result there will be convergence not to the equilibrium point, but to its 
neighborhood. Size of the neighborhood depends on $T$, $\gamma >0$ and the 
number of players.

\begin{figure}
\centering

\includegraphics[width=0.5\textwidth]{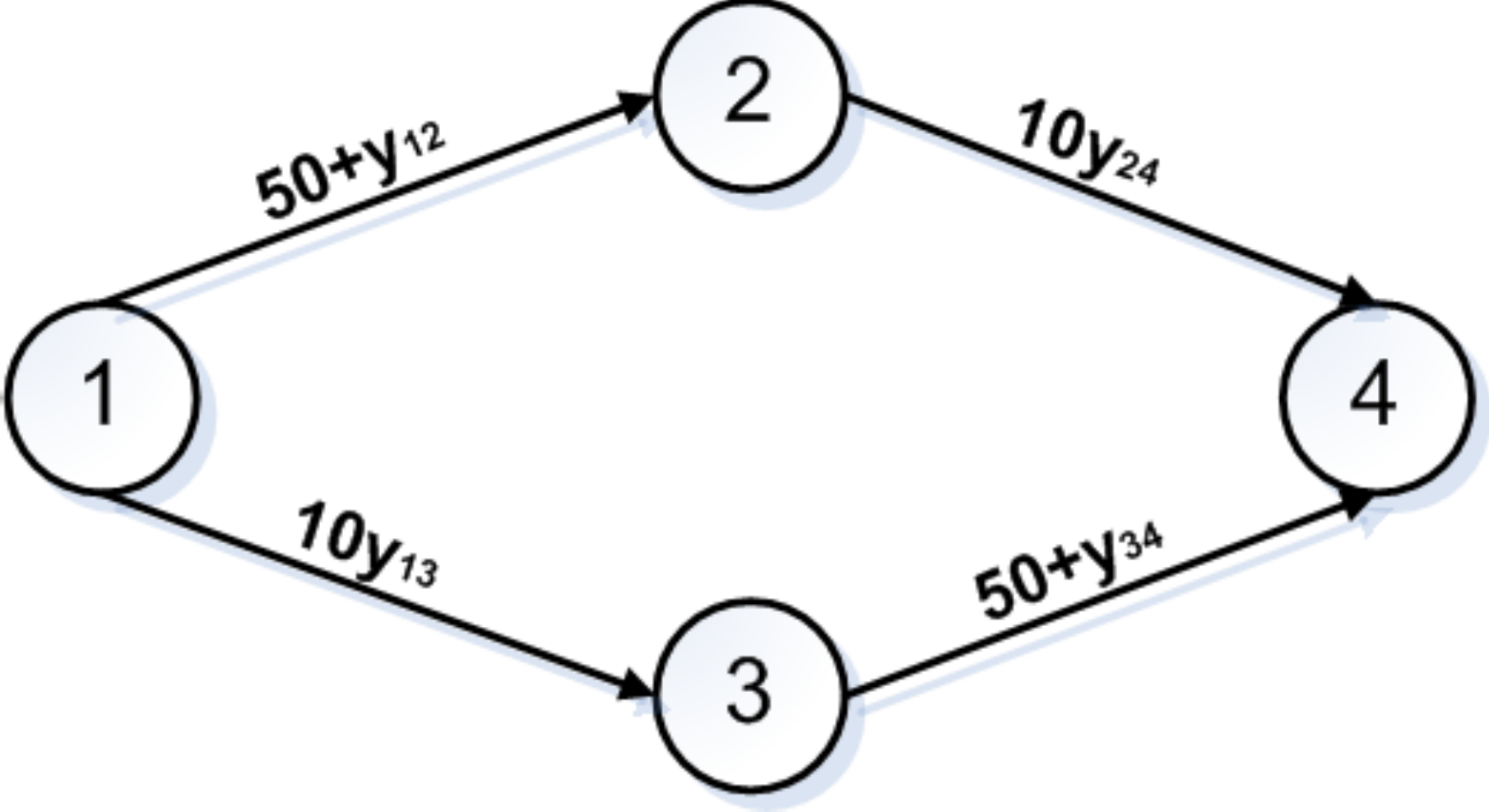}
\caption{\textbf{Case 1.} $x_{124} =x_{134} =3$. Total time for each path is $T=83$ min.}
\label{fig:1}       
\end{figure}

\begin{figure}
\centering

\includegraphics[width=0.5\textwidth]{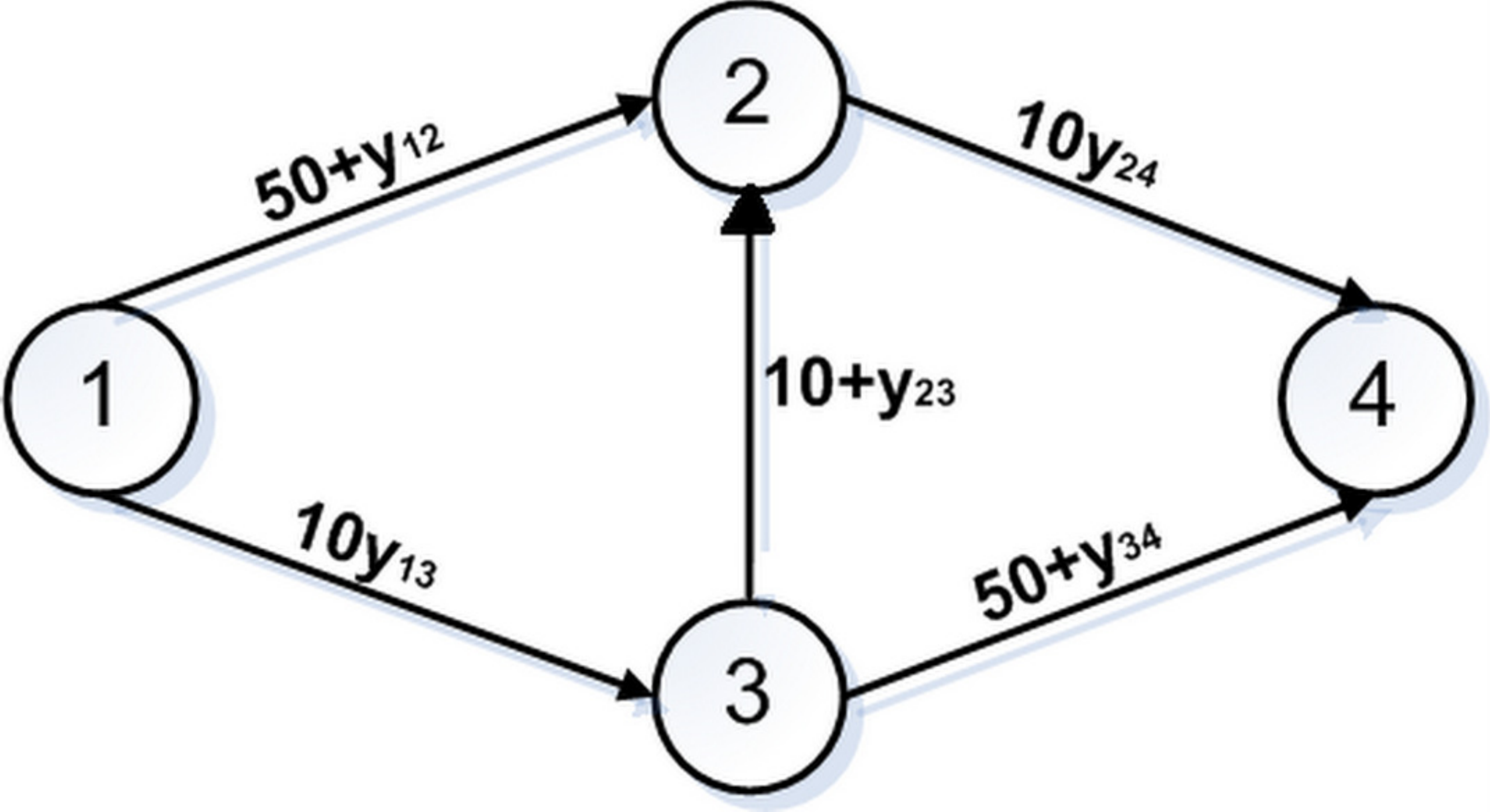}
\caption{\textbf{Case 2.} $x_{124} =x_{1234} =x_{134} =2$. Total time for each path is $T=92$ min.}
\label{fig:2}       
\end{figure}

\textbf{Example 1 (Braess paradox, 1968 \cite{15}) }
Let correspondence $x_{14} =6$ thousand cars/hour (see graphs on Figures \ref{fig:1} and \ref{fig:2}). Weight of the edges is time delay (in minutes) 
when the flow on the edge is $y_{ij} $ (thousand cars/hour). For example, in case 2 (see Figure \ref{fig:2}): $y_{24} =x_{1324} +x_{124} $. It is natural that time delay (at each 
of the edge) is a growth function of flow.

The following example shows, that under the very natural conditions 
vector-function of cost of the travel $\vec{G}\left( {\vec{x}} \right)$ 
can't be strictly monotone:
\[
\exists \vec{x},\vec{y}\in X\;\left( \vec{x}\ne \vec{y}\right)\colon\vec{G}\left( {\vec{x}} \right)=\vec{G}\left( {\vec{y}} 
\right)\Rightarrow \left\langle {\vec{G}\left( {\vec{x}} \right)-\;\vec{G}\left( {\vec{y}} \right),\;\vec{x}-\vec{y}} \right\rangle =0.
\]
This, for example, can be because of
\[
\vec{G}\left( {\vec{x}} \right)=\Theta ^T\vec{\tau }\left( {\vec{y}} 
\right),
\quad
\vec{y}=\Theta \vec{x},
\]
where $\vec{y}=\left\{ {y_e } \right\}_{e\in E} $ describes the loading of 
edges (arcs) of a graph of the transport network, $\vec{\tau }\left( {\vec{y}} \right)=\left\{ {\tau _e \left( {y_e } \right)} \right\}_{e\in E} $ -- 
vector-function of cost of the travel on the edges of transport network, 
$\Theta $ -- incidence matrix of edges and paths, and different vectors of 
flow distributions $\vec{x}$ may correspond to the same vector $\vec{y}=\Theta \vec{x}$.

\begin{figure}
\centering
\includegraphics[width=0.5\textwidth]{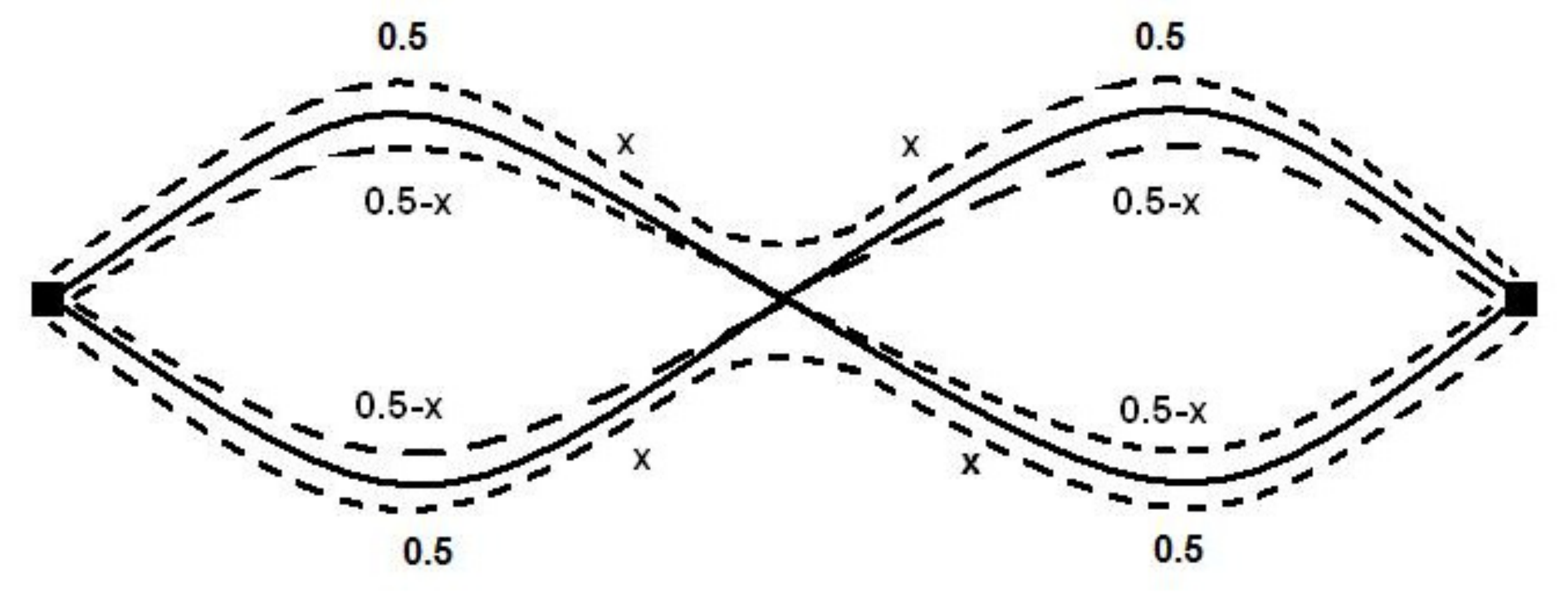}
\caption{Nonuniqueness of the equilibrium.}
\label{fig:3}
\end{figure}

\textbf{Example 2 (Nonuniqueness of the equilibrium; Shvetsov, 2010).} \\
On Figure \ref{fig:3} the equilibrium flow distribution is shown, for all $x\in \left[ 
{0,\;0.5} \right]$.

\begin{theorem}
\label{th:2}
Let $T>0$ be sufficiently small, $\sum\limits_{n=1}^\infty {\gamma _n } 
=\infty , \sum\limits_{n=1}^\infty {\left( {\gamma _n } \right)^2} <\infty 
$. Then
$\Psi \left( {\vec{x}\left( n \right)} 
\right)\mathrel{\mathop{\kern0pt\longrightarrow}\limits_{n\to \infty 
}^{a.s.}} \Psi \left( {\vec{x}^\ast } \right)$ and $\vec{x}\left( n 
\right)\mathrel{\mathop{\kern0pt\longrightarrow}\limits_{n\to \infty 
}^{a.s.}} \vec{x}^\ast \left( {\vec{x}\left( 0 \right)} \right).$
Note most of the elements of $\vec{x}^\ast \left( {\vec{x}\left( 0 \right)} \right)$ can be equal to zero.
\end{theorem}

We should notice, that Theorem \ref{th:2} is a refutation (in case of the considered dynamics) of the hypothesis \cite{16}. It states that in the case of 
non-unique Nash-Wardrop equilibrium, the equilibrium is more likely to 
realize, and it is a solution of the following linear-enthropy programming 
problem
\[
\sum\limits_{w\in W} {\sum\limits_{p\in P_w } {\left( {x_p \ln \left( {{x_p 
} \mathord{\left/ {\vphantom {{x_p } {\left| {P_w } \right|}}} \right. 
\kern-\nulldelimiterspace} {\left| {P_w } \right|}} \right)-x_p } \right)} } 
\to \mathop {\min }\limits_{\vec{x}\in X,\;\Theta \vec{x}=\vec{y}^\ast } 
,
\]
where $\vec{y}^\ast $ -- is the unique solution of $V\left( {\vec{y}} 
\right)=\sum\limits_{e\in E} {\int\limits_0^{y_e } {\tau _e \left( z 
\right)dz} } \to \min _{\vec{y}=\Theta \vec{x},\;\vec{x}\in X} $.
\section{Sketch of the proof of the Theorem \ref{th:1}}

\begin{lemma}
\label{lemma:1}
Let $f\left( w \right)=-T\ln w$, where $w\in \left( {0,\;1} \right)$. $\alpha _i >0$ --- are some randomly chosen, but fixed parameters; $w_i \in \left( {0,\;1}\right)$.
Let us consider functions
$$
F_0 \left( {\vec{w}} \right)=\frac{\sum\limits_i {\alpha _i w_i f\left( 
{w_i } \right)} }{\sum\limits_i {\alpha _i w_i } },
\text{ and }
F_1 \left( {\vec{w}} \right)=\frac{\sum\limits_i {\alpha _i f\left( {w_i } 
\right)} }{\sum\limits_i {\alpha _i } }.
$$
Then $F_0 \left( {\vec{w}} \right)\le F_1 \left( {\vec{w}} \right)$, and the equality is attained only when
$$
w_1 =w_2 =\ldots=w^\ast .$$
 
\end{lemma}
\begin{Proof}
 The proof is based on the consequent usage of the inequality between harmonic mean and geometric mean and then Cauchy inequality.
\end{Proof}

\begin{lemma}
\label{lemma:2}
For any $\vec{x}\left( n \right)\in X\colon\vec{x}\left( n \right)\ne \vec{x}^\ast $ holds the following inequality
\begin{multline}
\left\langle {\mathop{grad}\Psi \left( \vec{x}( n)\right),E\left[ {\left. \vec{x}(n+1)-\vec{x}\left( n \right) \right|\vec{x}\left( n \right)} \right]} \right\rangle =\\
=\left\langle {\vec{G}\left( \vec{x}( n) \right),E\left[{\left. {\vec{x}( n+1)-\vec{x}(n)} \right|\vec{x}\left( n \right)} \right]} \right\rangle <0. 
\end{multline}
\end{lemma}
\begin{Proof}
 Without restricting the generality, we can assume that $\forall \;p\in P\to x_p \left( n \right)\ge 1$. Then
$$
\left\langle {\vec{G}\left( {\vec{x}\left( n \right)} \right),\;E\left[ 
{\left. {\vec{x}\left( {n+1} \right)-\vec{x}\left( n \right)} \right|\vec{x}\left( n \right)} \right]} \right\rangle =
$$
$$
=\gamma_n \sum\limits_{w\in W} {d_w \left[ {\frac{\sum\limits_{p\in 
P_w } {x_p \left( n \right)G_p \left( {\vec{x}\left( n \right)} \right)\exp\left(-\cfrac{G_p(\vec{x}(n))}{T}\right)} }{\sum\limits_{p\in P_w } {x_p \left( 
n \right)\exp\left(-\cfrac{G_p(\vec{x}(n))}{T}\right)} }-\frac{\sum\limits_{p\in 
P_w } {G_p \left( {\vec{x}\left( n \right)} \right)x_p \left( n \right)} 
}{\sum\limits_{p\in P_w } {x_p \left( n \right)} }} \right]}.
$$
From Lemma \ref{lemma:1} it follows that
$$
\left\langle {\vec{G}\left( {\vec{x}\left( n \right)} \right),E\left[{\left. {\vec{x}\left( {n+1} \right)-\vec{x}\left( n \right)} \right|\vec{x}\left( n \right)} \right]} \right\rangle \le 0.
$$
And the equality can be attained only on the equilibrium vector $\vec{x}^\ast $, which can not be under considered hypothesis.
\end{Proof}
\begin{remark}
 Lemma \ref{lemma:2} can be more specified. At some neighborhood of the equilibrium there exists $l>0$, such that
$$
\left\langle {\vec{G}\left( {\vec{x}\left( n \right)} \right),\;E\left[ 
{\left. {\vec{x}\left( {n+1} \right)-\vec{x}\left( n \right)} \right|\vec{x}\left( n \right)} \right]} \right\rangle \le -l\gamma _n \cdot \left( 
{\Psi \left( {\vec{x}\left( n \right)} \right)-\Psi \left( {\vec{x}^\ast } 
\right)} \right).
$$
Also, with some reserves, we can change constraint $\vec{x}\left( n \right)\in X$ to $\vec{x}\left( n \right)\ge \vec{0}$.

\end{remark}

\begin{lemma}[\cite{25}, Chapter 2.2]
 Let \[
\sum\limits_{n=1}^\infty {\gamma _n } =\infty ,
\quad
\sum\limits_{n=1}^\infty {\left( {\gamma _n } \right)^2} <\infty .
\]
Then
$$
\Psi \left( {\vec{x}\left( n \right)} 
\right)\mathrel{\mathop{\kern0pt\longrightarrow}\limits_{n\to \infty 
}^{\text{a.s.}}} \Psi \left( {\vec{x}^\ast } \right)
$$
and if the equilibrium is unique, then also
$$
\vec{x}\left( n \right)\mathrel{\mathop{\kern0pt\longrightarrow}\limits_{n\to \infty }^{\text{a.s.}}} \vec{x}^\ast .
$$
\end{lemma}
\begin{Proof}
  Lemma \ref{lemma:2} and Theorem 1 from chapter 2.2 of \cite{25} allow us to consider only the situation, when $\vec{x}\left( n \right)$ is close to $\vec{x}^\ast $. 
Then by Taylor formula we have:
$$
E\left[ {\left. {\Psi \left( {\vec{x}\left( {n+1} \right)} \right)}\right|\vec{x}\left( n \right)} \right]=\Psi \left(\vec{x}( n)\right)+$$
$$
+\left\langle {\mbox{grad}\;\Psi \left( {\vec{x}\left( n \right)} 
\right),\;E\left[ {\left. {\vec{x}\left( {n+1} \right)-\vec{x}\left( n 
\right)} \right|\vec{x}\left( n \right)} \right]} \right\rangle +{\rm 
O}\left( {\left( {\gamma _n } \right)^2} \right).
$$

If we take mathematical expectations from both sides of this equality, we will get that there exists sufficiently large $C>0$, such that
$$
E\left( {\Psi \left( {\vec{x}\left( {n+1} \right)} \right)} \right)-\Psi 
\left( {\vec{x}^\ast } \right)\le \left( {1-l\gamma _n } \right)\cdot 
\left( {E\left( {\Psi \left( {\vec{x}\left( n \right)} \right)} 
\right)-\Psi \left( {\vec{x}^\ast } \right)} \right)+C\cdot \left( {\gamma 
_n } \right)^2.
$$
From more general statement from \cite{25} we get that
$$
E\left( {\Psi \left( {\vec{x}\left( n \right)} \right)} \right)-\Psi \left( 
{\vec{x}^\ast } 
\right)\mathrel{\mathop{\kern0pt\longrightarrow}\limits_{n\to \infty }} 0,
$$
if
$$
\sum\limits_{n=1}^\infty {\gamma _n } =\infty ,\quad\sum\limits_{n=1}^\infty {\left( {\gamma _n } \right)^2} <\infty .
$$
From Kolmogorov inequality follows
$$
P\left( {\forall \;n\ge n_0 \to \Psi \left( {\vec{x}\left( n \right)} 
\right)-\Psi \left( {\vec{x}^\ast } \right)\le \varepsilon } \right)\ge 
$$
$$
\ge 1-\varepsilon ^{-1}\cdot \left( {E\left( {\Psi \left( {\vec{x}\left( 
{n_0 } \right)} \right)-\Psi \left( {\vec{x}^\ast } \right)} 
\right)+\sum\limits_{k=n_0 }^\infty {\left( {\gamma _k } \right)^2} } 
\right).
$$
Which concludes the proof.
\end{Proof}
In the end we will formulate a known result, which is in high correlation with the proved one.

\begin{theorem}
Let\footnote{See formula (2.32), statement 2.2 and example on the page 1586 in \cite{26}} $T>0$. Then there $\exists\, C,\alpha >0\colon\forall \;N\in {\rm N}$ 
$$P\left( {\Psi \left({\frac{1}{N}\sum\limits_{n=1}^N {\vec{x}\left( n \right)} } \right)-\Psi 
_{\min } \ge \frac{\Omega }{\sqrt N }} \right)\le 2\exp \left( {-C\cdot 
\Omega } \right),$$ where $\gamma =\frac{\alpha }{\sqrt N }$

\end{theorem}

Authors thank S.A. Avvakumov, I.B. Gnedkov, Y.V. Dorn, I.S. Menshikov, E.A. 
Nurminskiy, A.A. Shananin, V.I. Shvetsov and especially A.V. Gasnikov.

The work was supported by RFBR 10-01-00321-a, 11-01-00494-a 
11-07-00162-a. The second author is partially supported by the Laboratory for Structural Methods of Data Analysis in Predictive Modeling, MIPT, RF government grant, ag. 11.G34.31.0073
%
%

%
%


\end{document}